\documentclass[12pt]{article}%
\usepackage{amsmath}
\usepackage{graphicx}%
\usepackage{amsfonts}%
\usepackage{amssymb}
\newtheorem{theorem}{Theorem}

\newtheorem{corollary}[theorem]{Corollary}

\newtheorem{proposition}[theorem]{Proposition}
\newtheorem{remark}[theorem]{Remark}

\begin{document}
\bigskip

\bigskip

\bigskip

\bigskip

\bigskip

\bigskip

\bigskip

\bigskip

\bigskip

\bigskip

\bigskip

\bigskip

\bigskip

\bigskip

\bigskip

\bigskip

\bigskip

\bigskip

\bigskip\bigskip

\bigskip\bigskip

\bigskip

\bigskip

{\Huge MARKOV LOOPS, DETERMINANTS }

{\Huge AND GAUSSIAN\ FIELDS\ \ }

{\Huge \bigskip}

Yves Le Jan

Math\'{e}matiques

Universit\'{e} Paris 11.

91405 Orsay. France

\bigskip

\bigskip yves.lejan@math.u-psud.fr

\section{Introduction}

The purpose of this article is to explore some simple relations between loop
measures, spanning trees, determinants, and Gaussian Markov fields. These
relations are related to Dynkin's isomorphism (cf \cite{Dy}, \cite{MR},
\cite{LJ1}) . Their potential interest could be suggested by noting that loop
measures were defined in \cite{LW} for planar Brownian motion and are related
to SLE processes (see also \cite{WW}). It is also the case for the free field
as shown in \cite{SCSH}. We present the results in the elementary framework of
symmetric Markov chains on a finite space, and then indicate how they can be
extended to more general Markov processes such as the two dimensional Brownian motion.

\section{Symmetric Markov processes on finite spaces}

Notations: Functions on finite (or countable) spaces are often denoted as
vectors and measures as covectors in coordinates with respect to the canonical
bases associated with points (the dual base being given by Dirac measures
$\delta_{x}$).

The multiplication operators defined by a function, $f$\ acting on functions
or on measures are in general simply denoted by $f$, but sometimes
multiplication operators by a function $f$ or a measure $\lambda$will be
denoted $M_{f}$ or $M_{\lambda}$. The function obtained as the density of a
measure $\mu$ with respect to some other measure $\nu$ is simply denoted
$\frac{\mu}{\nu}$.

\subsection{Energy and Markovian semigroups}

Let us first consider for simplicity the case of a symmetric irreducible
Markov chain with exponential holding times on a finite space $X$, with
generator $L_{y}^{x}=q^{x}(P_{y}^{x}-\delta_{y}^{x})$, $\lambda_{x},x\in X$
being a positive measure and $P$ a $\lambda$-symmetric stochastic transition
matrix: $\lambda_{x}P_{y}^{x}=\lambda_{y}P_{x}^{y}$ with $P_{x}^{x}=0$\ for
all $x$ in $X$.

We denote $P_{t}$ the semigroup $\exp(Lt)=\sum\frac{t^{k}}{k!}L^{k}$ and by
$m_{x}$ the measure $\frac{\lambda_{x}}{q^{x}}$. $L$ and $P_{t}$ are $m$-symmetric.

Recall that for any complex function $z^{x},x\in X$, the ''energy''%
\[
e(z)=\left\langle -Lz,\overline{z}\right\rangle _{m}=\sum_{x\in X}%
-(Lz)^{x}\overline{z}^{x}m_{x}%
\]
is nonnegative as it can be written
\[
e(z)=\frac{1}{2}\sum_{x,y}C_{x,y}(z^{x}-z^{y})(\overline{z}^{x}-\overline
{z}^{y})+\sum_{x}\kappa_{x}z^{x}\overline{z}^{x}=\sum_{x}\lambda_{x}%
z^{x}\overline{z}^{x}-\sum_{x,y}C_{x,y}z^{x}\overline{z}^{y}%
\]
with $C_{x,y}=C_{y,x}=\lambda_{x}P_{y}^{x}$ and $\kappa_{x}=\lambda_{x}%
(1-\sum_{y}P_{y}^{x})$, i.e. $\lambda_{x}=\kappa_{x}+\sum_{y}C_{x,y}%
=e(1_{\{x\}})$.

We say $\ (x,y)$ is a link iff $C_{x,y}>0$. An important exemple is the case
of a graph: Conductances are equal to zero or one and the conductance matrix
is the incidence matrix of the graph. 

\bigskip

\ The (complex) Dirichlet space $\mathbb{H}$\ is the space of
complex\ functions equipped with the energy scalar product defined by
polarisation of $e$. Note that the non negative symmetric ''conductance
matrix'' $C$ and the non negative equilibrium or ''killing'' measure $\kappa$
are the free parameters of the model. (so is $q$ but we will see it is
irrelevant for our purpose and we will mostly take it equal to $1$). The
lowest eigenvector of $-L$\ is nonnegative by the well known argument which
shows that the modulus contraction $z\rightarrow\left|  z\right|  $ lowers the
energy. We will assume (although it is not always necessary) the corresponding
eigenvalue is positive which means there is a ''mass gap'': For some positive
$\varepsilon$, the energy $e(z)$ dominates $\varepsilon\left\langle
z,\overline{z}\right\rangle _{m}$ for all $z$.

We denote by $V$ the associated potential operator $(-L)^{-1}=\int_{0}%
^{\infty}P_{t}dt$. They can be expressed in terms of the spectral resolution
of $L$.

We denote by $G$ the Green function defined on $X^{2}$ as $G^{x,y}%
=\frac{V_{y}^{x}}{m_{y}}=\frac{1}{\lambda_{y}}[(I-P)^{-1}]_{y}^{x}$ i.e.
$G=(M_{\lambda}-C)^{-1}$. It verifies $e(f,G\mu)=\left\langle f,\mu
\right\rangle $ for all function $f$ and measure $\mu$. In particular
$G\kappa=1$.

Different Markov chains associated to the same energy are equivalent under
time change. If $g$ is a positive function on $X$, in the new time scale
$\int_{0}^{t}g_{\xi_{s}}ds$, we obtain a Markov chain with $gm$%
-symmetric\ generator $\frac{1}{g}L$. Objects invariant under time change are
called \textsl{intrinsic}. The energy $e$, $P$ and the Green function $G$ are
obviously intrinsic but $L$, $V$ and $P_{t}$ are not. We will be interested
only in intrinsic objects. In this elementary framework, it is possible to
define a natural canonical time scale by taking $q=1$, but it will not be true
on continuous spaces.

\subsection{Recurrent chain}

Assume for simplicity that $q=1$. It will be convenient to add a cemetery
point $\Delta$ to $X$, and extend $C$, $\lambda$ and $G$ to $X^{\Delta
}=\{X\cup\Delta\}$by setting $C_{x,\Delta}=\kappa_{x}$ , $\lambda_{\Delta
}=\sum_{x\in X}\kappa_{x}$.\ and $G^{x,\Delta}=0$. Note that $\lambda
(X^{\Delta})=\sum_{X\times X}C_{x,y}+2\sum_{X}\kappa_{x})$

\bigskip

One can consider the recurrent ''resurrected'' Markov chain defined by the
extensions the conductances to $X^{\Delta}$. An energy $e^{R}$\ is defined by
the formula
\[
e^{R}(z)=\frac{1}{2}\sum_{x,y}C_{x,y}(z^{x}-z^{y})(\overline{z}^{x}%
-\overline{z}^{y})
\]
We denote by $P^{R}$ the transition kernel on $X^{\Delta}$\ defined by%

\[
e^{R}(z)=\left\langle z-P^{R}z,\overline{z}\right\rangle _{\lambda}%
\]
or equivalently by%
\[
\lbrack P^{R}]_{y}^{x}=\frac{C_{x,y}}{\sum_{y\in X^{\Delta}}C_{x,y}%
}=\frac{C_{x,y}}{\lambda_{x}}%
\]
Note that $P^{R}1=1$ so that $\lambda$ is now an invariant measure. Let
$\lambda^{\perp}$ be the space of functions on $X^{\Delta}$\ of zero $\lambda
$\ measure and by $V^{R}$\ the inverse of the restriction of $I-P^{R}$ to
$\lambda^{\perp}$.It vanishes on constants and has a mass gap on
$\lambda^{\perp}$. Setting for any signed measure $\nu$ of total charge zero
$G^{R}\nu=V^{R}\frac{\nu}{\lambda}$. we have for any function $f$,
$\left\langle \nu,f\right\rangle =e^{R}(G^{R}\nu,f)$ and in particular$\ f^{x}%
-f^{y}=e^{R}(G^{R}(\delta_{x}-\delta_{y}),f)$.

Note that for $\mu\in\lambda^{\perp}$ and carried by $X$, for all $x\in X$,
$\mu_{x}=e^{R}(G^{R}\mu,1_{x})=\lambda_{x}((I-P)G^{R}\mu)(x)-\kappa_{x}%
G^{R}\mu(\Delta)$. Hence, applying $G$ , it follows that on $X$,$\ G^{R}%
\mu=G^{R}\mu(\Delta)G\kappa+G\mu=G^{R}\mu(\Delta)+G\mu$. Moreover, as
$G^{R}\mu$ is in $\lambda^{\perp}$, $G^{R}\mu(\Delta)\lambda(X^{\Delta}%
)+\sum_{x\in X}\lambda_{x}(G\mu)_{x}=0$.

Therefore, $G^{R}\mu(\Delta)=\frac{-\left\langle \lambda,G\mu\right\rangle
}{\lambda(X^{\Delta})}$ and $G^{R}\mu=\frac{-\left\langle \lambda
,G\mu\right\rangle }{\lambda(X^{\Delta})}+G\mu$

\subsection{\bigskip Transfer matrix}

We can define a scalar product on the space $\mathbb{A}$ of antisymmetric
functions on $X^{\Delta}\times X^{\Delta}$\ as follows

$\left\langle \omega,\eta\right\rangle =\sum_{x,y}C_{x,y}\omega^{x,y}%
\eta^{x,y}.$ Denoting as in \cite{Lyo2}\ \ $df^{u,v}=f^{u}-f^{v}$, we note
that $\left\langle df,dg\right\rangle =e^{R}(f,g)$ In particular
\[
\left\langle df,dG^{R}(\delta_{x}-\delta_{y})\right\rangle =df^{x,y}%
\]
As the antisymmetric functions $df$ span the space of antisymmetric functions,
it follows that the scalar product is positive definite.

The symmetric transfer matrix $K$, indexed by pairs of oriented links, is
defined to be
\[
K^{(x,y),(u,v)}=G^{R}(\delta_{x}-\delta_{y})^{u}-G^{R}(\delta_{x}-\delta
_{y})^{v}=<dG^{R}(\delta^{x}-\delta^{y}),dG^{R}(\delta^{u}-\delta^{v})>
\]
for $x,y,u,v\in X^{\Delta}$, with $x\neq y,u\neq v$.

We see that for $x$\ and $y$ in $X$, $G^{R}(\delta_{x}-\delta_{y})^{u}%
-G^{R}(\delta_{x}-\delta_{y})^{v}=G(\delta_{x}-\delta_{y})^{u}-G(\delta
_{x}-\delta_{y})^{v}$.

We can see also that $G^{R}(\delta_{x}-\delta_{\Delta})=G\delta_{x}%
-\frac{-\left\langle \lambda,G\delta_{x}\right\rangle }{\lambda(X^{\Delta})}$.
So the same identity holds in $X^{\Delta}$.

Therefore, as $G^{x,\Delta}=0$, in all cases,
\[
K^{(x,y),(u,v)}=G^{x,u}+G^{y,v}-G^{x,v}-G^{y,u}%
\]

For every oriented link $\xi=(x,y)$ in \ $X^{\Delta}$,set $K^{\xi}%
=dG^{R}(\delta^{x}-\delta^{y})=dG(\delta^{x}-\delta^{y})$.

We have $\left\langle K^{\xi},K^{\eta}\right\rangle =K^{\xi,\eta}.$ $K$ will
be viewed as a linear operator on $\mathbb{A}$, self adjoint with respect to
$\left\langle \cdot,\cdot\right\rangle .$ (It can also be viewed as symmetric
with respect to the euclidean scalar product if we wish to use it Then it
appears as the inverse of the operator defined by $\left\langle \cdot
,\cdot\right\rangle $).

\section{Loop measures}

\subsection{Definitions}

For any integer $k$, let us define a based loop with $p$ points in $X$ as a
couple $(\xi,\tau)=((\xi_{m},1\leq m\leq p),(\tau_{m},1\leq m\leq
p+1),\mathbb{)}$ in $X^{p}\times\mathbb{R}_{+}^{p+1}$, and set $\xi_{1}%
=\xi_{p+1}$. $p$ will be denoted $p(\xi)$.

Based loops have a natural time parametrisation $\xi(t)$ and a time period
$T(\xi)=\sum_{i=1}^{p(\xi)+1}\tau_{i}$. If we denote $\sum_{i=1}^{m}\tau_{i}%
$\ by $T_{m}$: $\xi(t)=\xi_{m-1}$ on $[T_{m-1},T_{m})$ (with by
convention\ $T_{0}=0$ and $\xi_{0}=\xi_{p}$).

A $\sigma$-finite measure $\mu_{0}$ is defined on based loops by%
\[
\mu_{0}=\sum_{x\in X}\int_{0}^{\infty}\frac{1}{t}\mathbb{P}_{t}^{x,x}dt
\]

where $\mathbb{P}_{t}^{x,x}$ denotes the (non normalized) ''law'' of a path
from $x$ to $x$\ of duration $t$ : If $\sum_{i=1}^{h+1}t_{i}=t$,%
\[
\mathbb{P}_{t}^{x,x}(\xi(t_{1})=x_{1},...,\xi(t_{h})=x_{h})=[P_{t_{1}}%
]_{x_{1}}^{x}[P_{t_{2}-t_{1}}]_{x_{2}}^{x_{1}}...[P_{t-t_{h}}]_{x}^{x_{h}}%
\]
Note also that%
\begin{align*}
\mathbb{P}_{t}^{x,x}(p  &  =k,\xi_{2}=x_{2},...,\xi_{_{k}}=x_{k},T_{1}\in
dt_{1},...,T_{k}\in dt_{k})\\
&  =[P]_{x_{2}}^{x}[P]_{x_{3}}^{x_{2}}...[P]_{x}^{x_{k}}1_{\{0<t_{1}%
<...t_{k}<t\}}q^{x}e^{-q^{x}t_{1}}...q_{x_{k}}e^{-q_{x_{k}}(t_{k}-t_{k-1}%
)}e^{-q_{x}(t-t_{k})}dt_{1}...dt_{k}%
\end{align*}
A loop is defined as an equivalence class of based loops for the $\mathbb{R}%
$-shift that acts naturally. $\mu_{0}$ is shift invariant, It induces
a\textsl{ measure }$\mu$\textsl{\ on loops}.

Note also that the measure $d\widetilde{\mu}_{0}=\frac{Tq_{\xi_{1}}}{\int
_{0}^{T}q_{\xi(s)}ds}d\mu_{0}$ which is not shift invariant also induces $\mu$
on loops.

It writes%
\begin{align*}
\widetilde{\mu}_{0}(p(\xi)  &  =k,\xi_{1}=x_{1},...,\xi_{k}=x_{k},T_{1}\in
dt_{1},...,T_{k}\in dt_{k},T\in dt)\\
&  =[P]_{x_{2}}^{x_{1}}[P]_{x_{3}}^{x_{2}}...[P]_{x}^{x_{k}}\frac{1_{\{0<t_{1}%
<...<t_{k}<t\}}}{\int_{0}^{t}q_{\xi(s)}ds}e^{-q^{x_{1}}t_{1}}e^{-q_{x_{2}%
}(t_{2}-t_{1})}...e^{-q_{x}(t-t_{k})}q_{x_{1}}dt_{1}...q_{x_{k}}dt_{k}%
q_{x_{1}}dt
\end{align*}
for $k\geq2$ and%
\[
\widetilde{\mu}_{0}\{p(\xi)=1,\xi_{1}=x,\tau_{1}\in dt_{1}\}=\frac{e^{-q_{x}%
t_{1}}}{t_{1}}dt_{1}%
\]

It is clear, in that form, that a time change transforms the $\widetilde{\mu
}_{0}$'s of Markov chains associated with the same energy one into each other,
and therefore the same holds for $\mu$: this is analogous to conformal
invariance. Hence the restriction $\mu_{I}$\ of $\mu$ to the $\sigma$-field of
sets of loops invariant by time change (i.e. intrinsic sets) is intrinsic. It
depends only on $e$. As we are interested in the restriction $\mu_{I}$ of
$\mu$ to intrinsic sets, from now on we will denote simply $\mu_{I}$ by $\mu$

Intrinsic sets are defined by the discrete loop $\xi_{i}$ (in circular order,
up to translation) and the associated intrinsic times $\frac{\tau_{i}}{m_{i}%
}=\tau_{i}^{\ast}$. Conditionally to the discrete loop, these are independent
exponential variables with parameters $\lambda_{i}$.%

\begin{equation}
\mu=\sum_{x\in X}e^{-\lambda_{x}\tau^{\ast}}\frac{d\tau^{\ast}}{\tau^{\ast}%
}+\sum_{p=2}^{\infty}\sum_{(\xi_{i},i\in\mathbb{Z}/p\mathbb{Z)\in}X^{p}}%
\prod_{i\in\mathbb{Z}/p\mathbb{Z}}C_{\xi_{i},\xi_{i+1}}e^{-\lambda_{\xi_{i}%
}\tau_{i}^{\ast}}d\tau_{i}^{\ast} \label{Def}%
\end{equation}

\bigskip

Sets of\ discrete loop are the most important intrinsic sets, though we will
see that to establish a connection with Gaussian fields it is important to
consider occupation times. The simplest intrinsic variables are
\[
N_{x,y}=\#\{i:\xi_{i}=x,\xi_{i+1}=y\}
\]
and%
\[
N_{x}=\sum_{y}N_{x,y}%
\]
Note that $N_{x}=\#\{i\geq1:\xi_{i}=x\}$\ except for trivial one point loops.

A \textsl{bridge measure }$\mu^{x,y}$ can be defined on paths $\gamma$ from
$x$ to $y$: $\mu^{x,y}(d\gamma)=\frac{1}{m_{y}}\int_{0}^{\infty}\mathbb{P}%
_{t}^{x,y}(d\gamma)dt$ with%

\[
\mathbb{P}_{t}^{x,y}(\gamma(t_{1})=x_{1},...,\gamma(t_{h})=x_{h})=P_{t_{1}%
}(x,x_{1})P_{t_{2}-t_{1}}(x_{1},x_{2})...P_{t-t_{h}}(x_{h},y)
\]
Note that the mass of $\mu^{x,y}$ is $\frac{V_{y}^{x}}{m_{y}}=G^{x,y}$. We
also have, with similar notations as the one defined for loops%

\begin{align*}
\mu^{x,y}(p(\gamma)  &  =k,\gamma_{2}=x_{2},...,\gamma_{k-1}=x_{k-1},T_{1}\in
dt_{1},...,T_{k-1}\in dt_{k-1},T\in dt)\\
&  =\frac{C_{x,x_{2}}C_{x_{2},x_{3}}...C_{x_{k-1},y}}{\lambda_{x}%
\lambda_{x_{2}}...\lambda_{y}}1_{\{0<t_{1}<...<t_{k}<t\}}e^{-q_{x}t_{1}%
}e^{-q_{x_{2}}(t_{2}-t_{1})}...e^{-q_{y}(t-t_{k})}q_{x}dt_{1}...q_{x_{k-1}%
}dt_{k}q_{y}dt
\end{align*}

so that the restriction of $\mu^{x,y}$ to intrinsic sets of paths is intrinsic.

Finally, we denote $\mathbb{P}^{x}$\ the family of probability laws on paths
defined by $P_{t}$.%
\[
\mathbb{P}^{x}(\gamma(t_{1})=x_{1},...,\gamma(t_{h})=x_{h})=P_{t_{1}}%
(x,x_{1})P_{t_{2}-t_{1}}(x_{1},x_{2})...P_{t_{h}-t_{h-1}}(x_{h-1},x_{h})
\]%

\begin{align*}
\mathbb{P}_{x}(p(\gamma)  &  =k,\gamma_{2}=x_{2},...,\gamma_{k}=x_{k},T_{1}\in
dt_{1},...,T_{k}\in dt_{k})\\
&  =\frac{C_{x,x_{2}}...C_{x_{k-1},x_{k}}\kappa_{x_{k}}}{\lambda_{x}%
\lambda_{x_{2}}...\lambda_{x_{k}}}1_{\{0<t_{1}<...<t_{k}\}}e^{-q_{x}t_{1}%
}...e^{-q_{x_{k}}(t_{k}-t_{k-1})}q_{x}dt_{1}...q_{x_{k}}dt_{k}%
\end{align*}

\subsection{First properties}

If $D$ is a subset of $X$, the restriction of $\mu$ to loops contained in $D$,
denoted $\mu^{D}$ is clearly the loop measure induced by the Markov chain
killed at the exit of $D$. This can be called the \textsl{restriction
property}.

Let us recall that this killed Markov chain is defined by the restriction of
$\lambda$ to $D$ and the restriction $P^{D}$ of $P\ $\ to $D^{2}\ $(or
equivalently by the restriction $e_{D}$ of the Dirichlet norm $e$ to functions
vanishing outside $D$) and (for the time scale), by the restriction of $q$ to
$D$.

\bigskip

From now on in this section, we will take $q_{x}=1$ for all $x$. Then $\mu
_{0}$ takes a simpler form:%

\begin{align*}
\mu_{0}(p(\xi)  &  =k,\xi_{1}=x_{1},...,\xi_{k}=x_{k},T_{1}\in dt_{1}%
,...,T_{k}\in dt_{k},T\in dt)\\
&  =P_{x_{2}}^{x_{1}}...P_{x_{1}}^{x_{k}}\frac{1_{\{0<t_{1}<...<t_{k}<t\}}}%
{t}e^{-t}dt_{1}...dt_{k}dt
\end{align*}
for $k>1$ and $\mu_{0}\{p(\xi)=1,\xi_{1}=x_{1},\tau_{1}\in dt_{1}%
\}=\frac{e^{-t_{1}}}{t_{1}}dt_{1}$

It follows that for $k>0$,
\[
\mu_{0}(p(\xi)=k,\xi_{1}=x_{1},...,\xi_{k}=x_{k})=\frac{1}{k}P_{x_{2}}^{x_{1}%
}...P_{x_{1}}^{x_{k}}=\frac{1}{k}\prod_{x,y}C_{x,y}^{N_{x,y}}\prod_{x}%
\lambda_{x}^{-N_{x}}%
\]
as$\int\frac{t^{k-1}}{k!}e^{-t}dt=\frac{1}{k}$ and conditionally to
$p(\xi)=k,\xi_{1}=x_{1},...,\xi_{k}=x_{k}$, $T$ is a gamma variable of density
$\frac{t^{k-1}}{(k-1)!}e^{-t}$ on $\mathbb{R}_{+}$\ and $(\frac{T_{i}}{T}1\leq
i\leq k)$ an independent ordered $k-$sample of the uniform distribution on
$(0,1)$.

In particular, we obtain that, for $k\geq2$
\[
\mu(p=k)=\mu_{0}(p=k)=\frac{1}{k}Tr(P^{k})
\]
and therefore, as $Tr(P)=0$,
\[
\mu(p>0)=-\log(\det(I-P))=-\log(\frac{\det(G)}{\prod_{x}\lambda_{x}})
\]

as denoting $M_{\lambda}$ the diagonal matrix with entries $\lambda_{x}$,
$\det(I-P)=\frac{\det(M_{\lambda}-C)}{\det(M_{\lambda})}$. Moreover
\[
\int p(l)\mu(dl)=Tr((I-P)^{-1}P)
\]

Similarly, for any $x\neq y$ in $X$ and $s\in\lbrack0,1]$, setting
$P_{u,v}^{(s)}=P_{v}^{u}$ if $(u,v)\neq(x,y)$ and $P_{x,y}^{(s)}=sP_{y}^{x}$,
we have:%

\[
\mu(s^{N_{x,y}}1_{\{p>0\}})=-\log(\det(I-P^{(s)}))
\]

Differentiating in $s=1$, it comes that%
\[
\mu(N_{x,y})=[(I-P)^{-1}]_{x}^{y}P_{y}^{x}=G^{x,y}C_{x,y}%
\]

and $\mu(N_{x})=\sum_{y}\mu(N_{x,y})=\lambda_{x}G^{x,x}-1$ (as $(M_{\lambda
}-C)G=Id$).

\section{Poisson process of loops and occupation field\label{occup}}

\subsection{Occupation field}

To each loop $l$ we associate an occupation field $\{\widehat{l_{x}},x\in X\}$
defined by
\[
\widehat{l}^{x}=\int_{0}^{T(l)}1_{\{\xi(s)=x\}}\frac{q^{\xi_{s}}}{m_{\xi(s)}%
}ds=\sum_{i=1}^{p(l)}1_{\{\xi_{i-1}=x\}}\frac{q^{x}\tau_{i}}{m_{x}}=\sum
_{i=1}^{p(l)}1_{\{\xi_{i-1}=x\}}\tau_{i}^{\ast}%
\]
for any representative $(\xi,\tau)$ of $l$. It is independent of the time
scale (i.e.''intrinsic'').

For a path $\gamma$, $\widehat{\gamma}$ is defined in the same way.

From now on\ \emph{we will take }$q=1$\emph{.}

Note that
\begin{equation}
\mu((1-e^{-\alpha\widehat{l}^{x}})1_{\{p=1\}})=\int_{0}^{\infty}%
(e^{-(\frac{\alpha}{\lambda_{x}}+1)t}-e^{-t})\frac{dt}{t}=\log(\frac{\lambda
_{x}}{\alpha+\lambda_{x}}) \label{points}%
\end{equation}

In particular, $\mu(\widehat{l}^{x}1_{\{p=1\}})=\frac{1}{\lambda_{x}}$.

From formula \ref{Def}, we get easily that for any function $\Phi$ of the
discrete loop and $k\geq1$,
\[
\mu((\widehat{l}^{x})^{k}1_{\{p>1\}}\Phi)=\mu((N_{x}+k-1)...(N_{x}+1)N_{x}%
\Phi)
\]

In particular, $\mu(\widehat{l}^{x})=\frac{1}{\lambda_{x}}[\mu(N_{x}%
)+1]=G^{x,x}$.

Note that functions of $\widehat{l}$ are not the only intrinsic functions.
Other intrinsic variables of interest are, for $k\geq2$

$\widehat{l}^{x_{1},...,x_{k}}=\frac{1}{k}\sum_{j=0}^{k-1}\int_{0<t_{1}%
<...<t_{k}<T}1_{\{\xi(t_{1})=x_{1+j},....\xi(t_{k-j})=x_{k},...\xi
(t_{k})=x_{j}\}}\prod\frac{1}{\lambda_{x_{i}}}dt_{i}$

$=\frac{1}{k}\sum_{j=0}^{k-1}\sum_{1\leq i_{1}<..<i_{k}\leq p(l)}\prod
_{l=1}^{k}1_{\{\xi_{i_{l}-1}=x_{l+j}\}}\tau_{i_{l}}^{\ast}$ and one can check
that $\mu(\widehat{l}^{x_{1},...,x_{k}})=G^{x_{1},x_{2}}G^{x_{2},x_{3}%
}...G^{x_{k},x_{1}}$. Note that in general $\widehat{l}^{x_{1},...,x_{k}}$
cannot be expressed in terms of $\widehat{l}$ for $k>3$.

For $x_{1}=x_{2}=...=x_{k}$, we obtain self intersection local times
$\widehat{l}^{x,k}=\sum_{1\leq i_{1}<..<i_{k}\leq p(l)}\prod_{l=1}^{k}%
1_{\{\xi_{i_{l}-1}=x\}}\tau_{i_{l}}^{\ast}$

For any function $\Phi$ of the discrete loop, $\mu(\widehat{l}^{x,2}\Phi
)=\mu(\frac{N_{x}(N_{x}-1)}{2}\Phi)$ since $\widehat{l}^{x,2}=\frac{1}%
{2}((\widehat{l}^{x})^{2}-\sum_{i=1}^{p(l)}1_{\{\xi_{i-1}=x\}}(\tau_{i}^{\ast
})^{2})$ and $\mu(\Phi\sum_{i=1}^{p(l)}1_{\{\xi_{i-1}=x\}}(\tau_{i}^{\ast
})^{2}))=2\mu(\Phi N_{x})$

More generally one proves in a similar way that $\mu(\widehat{l}^{x,k}%
\Phi)=\mu(\frac{N_{x}(N_{x}-1)...(N_{x}-k+1)}{k!}\Phi)$

\bigskip

From the Feynman-Kac formula, it comes easily that, denoting $M_{\frac{\chi
}{\lambda}}$ the diagonal matrix with coefficients $\frac{\chi_{x}}%
{\lambda_{x}}$

$\mathbb{P}_{x,x}^{t}(e^{-\left\langle \widehat{l},\chi\right\rangle }%
-1)=\exp(t(P-I-M_{_{\frac{\chi}{\lambda}}}))_{x,x}-\exp(t(P-I))_{x,x}$.
Integrating in $t$ after expanding, we get from the definition of $\mu$ (first
for $\chi$ small enough):%

\[
\int(e^{-\left\langle \widehat{l},\chi\right\rangle }-1)d\mu(l)=\sum
_{k=1}^{\infty}\frac{1}{k}[Tr((P-M_{_{\frac{\chi}{\lambda}}})^{k}%
)-Tr((P)^{k})]
\]

Hence%
\[
\int(e^{-\left\langle \widehat{l},\chi\right\rangle }-1)d\mu(l)=\log
[\det(-L(-L+M_{\chi/\lambda})^{-1})]=-\log\det(I+VM_{\frac{\chi}{\lambda}})
\]
which now holds for all non negative $\chi$. Set $V_{\chi}=(-L+M_{_{\frac{\chi
}{\lambda}}})^{-1}$ and $G_{\chi}=V_{\chi}M_{\frac{1}{\lambda}}$. It is an
intrinsic symmetric nonnegative function on $X\times X$. $G_{0}$ is the Green
function $G$, and $G_{\chi}$ can be viewed as the Green function of the energy
form $e_{\chi}=e+\left\|  {}\right\|  _{L^{2}(\chi)}^{2}$. Note that $e_{\chi
}$ has the same conductances $C$ as $e,$ but $\chi$ is added to the killing
measure. We have also the ''resolvent'' equation $V-V_{\chi}=VM_{\frac{\chi
}{\lambda}}V_{\chi}=V_{\chi}M_{\frac{\chi}{\lambda}}V$. Then, $G-G_{\chi
}=GM_{\chi}G_{\chi}=G_{\chi}M_{\chi}G$. Also:%

\begin{equation}
\det(I+GM_{\chi})^{-1}=\det(I-G_{\chi}M_{\chi})=\frac{\det(G_{\chi})}{\det(G)}
\label{F1}%
\end{equation}
Finally we have the

\begin{proposition}
\label{xx}i)$\mu(e^{-\left\langle \widehat{l},\chi\right\rangle }%
-1)=-\log(\det(I+GM_{\chi}))=\log(\det(I-G_{\chi}M_{\chi}))=\log(\det(G_{\chi
}G^{-1}))$
\end{proposition}

Note that in this calculation, the trace and the determinant are applied to
matrices indexed by $X$. Note also that $\det(I+GM_{\chi})=\det(I+M_{\sqrt
{\chi}}GM_{\sqrt{\chi}})$ and $\det(I-G_{\chi}M_{\chi})=\det(I-M_{\sqrt{\chi}%
}G_{\chi}M_{\sqrt{\chi}})$, so we can deal with symmetric matrices..

In view of generalizing them to continuous spaces in an intrinsic form (i.e.
in a form invariant under time change), , $G$ and $G_{\chi}$ will be
interpreted as symmetric elements of $\mathbb{H\otimes H}$, or as linear
operators from $\mathbb{H}^{\prime}$ into $\mathbb{H}$. $G$ is a canonical
bijection. $\frac{\det(G_{\chi})}{\det(G)}$ can be viewed as the determinant
of the operator $G_{\chi}G^{-1}$ acting on $\mathbb{H}$.

\subsection{Poisson process of loops}

Still following the idea of \cite{LW}, define, for all positive $\alpha$,\ the
Poisson process of loops $\mathcal{L}_{\alpha}$ with intensity $\alpha\mu$. We
denote by $\mathbb{P}$\ or $\mathbb{P}_{\mathcal{L}_{\alpha}}$ its
distribution. Note that by the restriction property, $\mathcal{L}_{\alpha}%
^{D}=\{l\in\mathcal{L}_{\alpha},l\subseteq D\}$\ is a Poisson process of loops
with intensity $\mu^{D}$, and that $\mathcal{L}_{\alpha}^{D}$ \ is independent
of $\mathcal{L}_{\alpha}\backslash\mathcal{L}_{\alpha}^{D}$.

We denote by $\mathcal{L}_{\alpha}^{d}$ the set of non trivial discrete loops
in $\mathcal{L}_{\alpha}$. Then, $\mathbb{P(}\mathcal{L}_{\alpha}^{d}%
=\{l_{1},l_{2},...l_{k}\})=e^{-\alpha\mu(p>0)}\alpha^{k}\frac{\mu(l_{1}%
)...\mu(l_{k})}{k!}=[\frac{\det(G)}{\prod_{x}\lambda_{x}}]^{\alpha}\prod
_{x,y}C_{x,y}^{N_{x,y}^{(\alpha)}}\prod_{x}\lambda_{x}^{-N_{x}^{(\alpha)}}$
with $N_{x}^{(\alpha)}=\sum_{l\in\mathcal{L}_{\alpha}}N_{x}(l)$ and
$N_{x,y}^{(\alpha)}=\sum_{l\in\mathcal{L}_{\alpha}}N_{x,y}(l)$.

\begin{remark}
\label{unif} It follows that the probability of a discrete loop configuration
depends only on the variables $N_{x,y}+N_{y,x}$, i.e. the total number of
traversals of non oriented links. In particular, it does not depend on the
orientation of the loops It should be noted that under loop or path measures,
the conditional distributions of discrete loops or paths given the values of
all $N_{x,y}+N_{y,x}$'s is uniform. The $N_{x,y}+N_{y,x}$ ($N_{x,y}$)
configuration can be called the associated random \ (oriented) graph. Note
however that any configuration of $N_{x,y}+N_{y,x}$ does not correspond to a
loop configuration.
\end{remark}

We can associate to $\mathcal{L}_{\alpha}$ the $\sigma$-finite measure
\[
\widehat{\mathcal{L}_{\alpha}}=\sum_{l\in\mathcal{L}_{\alpha}}\widehat{l}%
\]

Then, for any non-negative measure $\chi$ on $X$
\[
\mathbb{E}(e^{-\left\langle \widehat{\mathcal{L}_{\alpha}},\chi\right\rangle
})=\exp(\alpha\int(e^{-\left\langle \widehat{l},\chi\right\rangle }%
-1)d\mu(l))
\]
and%
\[
\mathbb{E}(e^{-\left\langle \widehat{\mathcal{L}_{\alpha}},\chi\right\rangle
})=[\det(-L(-L+M_{\chi/\lambda})^{-1})]^{\alpha}=\det(I+VM_{\frac{\chi
}{\lambda}})^{-\alpha}%
\]
Finally we have the

\begin{proposition}
\label{det1}$\mathbb{E}(e^{-\left\langle \widehat{\mathcal{L}_{\alpha}}%
,\chi\right\rangle })=\det(I+GM_{\chi})^{-\alpha}=\det(I-G_{\chi}M_{\chi
})^{\alpha}=\det(G_{\chi}G^{-1})^{\alpha}$
\end{proposition}

Many calculations follow from proposition \ref{xx}.

It follows that $\mathbb{E}(\widehat{\mathcal{L}_{\alpha}}^{x})=\alpha
G_{x,x}$ and we recover that $\mu(\widehat{l}_{x})=G_{x,x}$.

On loops and paths, we define the restricted intrinsic $\sigma$-field
$\mathcal{I}_{R}$\ as generated the variables $N_{x,y}$ with $y.$ possibly
equal to $\Delta$ in the case of paths, with $N_{x,\Delta}=0$ or $1$. from
(\ref{points}),\
\[
\mathbb{E}(e^{-\sum\chi_{i}\left\langle \widehat{\mathcal{L}_{\alpha}}%
,\delta_{x_{i}}\right\rangle }|\mathcal{I}_{R})=\prod_{i=1}^{k}(\frac{\lambda
_{x_{i}}}{\lambda_{x_{i}}+\chi_{i}})^{N_{x_{i}}^{(\alpha)}+1}%
\]

The distribution of $\{N_{x}^{(\alpha)},x\in X\}$ follows easily, in terms of
generating functions:
\[
\mathbb{E}(\prod_{i=1}^{k}(s_{i}^{N_{x_{i}}^{(\alpha)}+1})=\det(\delta
_{i,j}+\sqrt{\frac{\lambda_{x_{i}}\lambda_{j}(1-s_{i})(1-s_{j})}{s_{i}s_{j}}%
}G_{x_{i},x_{j}})^{-\alpha}%
\]

Note also that
\[
\mathbb{E}((\widehat{\mathcal{L}_{\alpha}}^{x})^{k}|\mathcal{I}_{R}%
)=\frac{(N_{x}^{(\alpha)}+k)(N_{x}^{(\alpha)}+k-1)...(N_{x}^{(\alpha)}%
+1)}{k!\lambda_{x}^{k}}%
\]
and if self intersection local times are defined as

$\widehat{\mathcal{L}_{\alpha}}^{x,k}=\sum_{m=1}^{k}\sum_{k_{1}+...+k_{m}%
=k}\sum_{l_{1}\neq l_{2}...\neq l_{m}\in\mathcal{L}_{\alpha}^{+}}\prod
_{j=1}^{m}\widehat{l_{j}}^{x,k_{j}}$, we get easily that
\[
\mathbb{E}(\widehat{\mathcal{L}_{\alpha}}^{x,k}|\mathcal{I}_{R})=\frac{1}%
{\lambda_{x}^{k}}(N_{x}^{(\alpha)}-k+1)...(N_{x}^{(\alpha)}-1)N_{x}^{(\alpha)}%
\]

\bigskip

Note also that since $G_{\chi}M_{\chi}$\ is a contraction, from determinant
expansions given in \cite{VJ1} and \cite{VJ2}, we have%

\[
\mathbb{E}(\left\langle \widehat{\mathcal{L}_{\alpha}},\chi\right\rangle
^{k})=\sum\chi_{i_{1}}...\chi_{i_{k}}Per_{\alpha}(G_{i_{l},i_{m}},1\leq
l,m\leq k)
\]

Here the $\alpha$-permanent $Per_{a}$ is defined as $\sum_{\sigma
\in\mathcal{S}_{k}}\alpha^{m(\sigma)}G_{i_{1},i_{\sigma(1)}}...G_{i_{k}%
,i_{\sigma(k)}}$ with $m(\sigma)$ denoting the number of cycles in \ $\sigma$.

Let $[H^{F}]_{\cdot}^{x}$ be the hitting distribution of $F$ by the Markov
chain starting at $F$. Set $D=F^{c}$ and denote $e^{D}$, $V^{D}%
=[(I-P)|_{D\times D}]^{-1}$ and $G^{D}=[(M_{\lambda}-C)|_{D\times D}]^{-1}$
the Dirichlet norm, the potential and the Green function of the process killed
at the hitting of $F$. Recall that $V=V^{D}+H^{F}V$ and $G=$ $G^{D}+H^{F}G$.

Taking $\chi=a1_{F}$ with $F$ finite, and letting $a$ increase to infinity, we get

$\lim_{a\uparrow\infty}(G_{\chi}M_{\chi})=H^{F}$ which is $I$ on $F$.
Therefore by proposition \ref{xx}, one checks that $\mathbb{P(}\widehat
{\mathcal{L}}_{\alpha}(F)=0)=\det(I-H^{F})=0$ and $\mu(\widehat{l(}%
F)>0)=\infty$. But this is clearly due to trivial loops as it can be seen
directly from the definition of $\mu$ that in this simple framework they cover
the whole space $X$.

Note however that $\mu(\widehat{l(}F)>0,p>0)=\mu(p>0)-\mu(\widehat{l(}F)=0,p>0)$

$=\mu(p>0)-\mu^{D}(p>0)=-\log(\frac{\det(I-P)}{\det_{D\times D}(I-P)}%
)=\log(\frac{\det(G^{D})}{\prod_{x\in F}\lambda_{x}\det(G)})$

It follows that the probability no non trivial loop (i.e.a loop which is not
reduced to a point) in $\mathcal{L}_{\alpha}$ intersects $F$ equals
$(\frac{\det(G^{D})}{\prod_{x\in F}\lambda_{x}\det(G)})^{\alpha}$

Recall that for any $(n+p,n+p)$ invertible matrix $A$, $\det(A^{-1}%
)\det(A_{ij}1\leq i,j\leq n)=\det(A^{-1})\det(Ae_{1},...Ae_{n},e_{n+1}%
,...e_{n+p})$

$=\det(e_{1},...e_{n},A^{-1}e_{n+1},...A^{-1}e_{n+p})=\det((A^{-1}%
)_{k,l},n\leq k,l\leq n+p)$.

In particular, $\det(G^{D})=\frac{\det(G)}{\det(G|_{F\times F})}$, so we have the

\begin{corollary}
The probability that no non trivial loop in $\mathcal{L}_{\alpha}$ intersects
$F$ equals $(\prod_{x\in F}\lambda_{x}\det_{F\times F}(G)^{-\alpha}$
\end{corollary}

In particular, it follows that the probability a non trivial loop in
$\mathcal{L}_{\alpha}$ visits $x$ equals $1-(\frac{1}{\lambda_{x}G^{x,x}%
})^{\alpha}$

Also, if $F_{1}$ and $F_{2}$ are disjoint, $\mu(\prod\widehat{l(}F_{i}%
)>0)=\mu(p>0)+\mu(\sum\widehat{l(}F_{i})=0,p>0)-\mu(\widehat{l(}%
F_{1})=0,p>0)-\mu(\widehat{l(}F_{2})=0,p>0)$

$=\log(\frac{\det(G)\det(G^{D_{1}\cap D_{2}})}{\det(G^{D_{1}})\det(G^{D_{2}}%
)})$ and this formula is easily generalized to $n$ disjoint sets.
\[
\mu(\prod\widehat{l(}F_{i})>0)=\log(\frac{\det(G)\prod_{i<j}\det(G^{D_{i}\cap
D_{j}})...}{\prod\det(G^{D_{i}})\prod_{i<j<k}\det(G^{D_{i}\cap D_{j}\cap
D_{k}})...}%
\]
The positivity yields an interesting determinant product inequality.

It follows in particular that the probability a non trivial loop in
$\mathcal{L}_{\alpha}$\ visits two distinct points $x$ and $y$ equals
$1-(\frac{G^{x,x}G^{y,y}-(G^{x,y})^{2}}{G^{x,x}G^{y,y}})^{\alpha}$ and
$\frac{(G^{x,y})^{2}}{G^{x,x}G^{y,y}}$\ if $\alpha=1$.

Note finally that if $\chi$ has support in $D$, by the restriction property%

\[
\mu(1_{\{\widehat{l(}F)=0\}}(e^{-<\widehat{l},\chi>}-1))=-\log(\det
(I+G^{D}M_{\chi}))=\log(\det(G_{\chi}^{D})[G^{D}]^{-1})
\]

Here the determinants are taken on matrices indexed by $D$. or equivalently on
operators on $\mathbb{H}^{D}$.

For paths we have $\mathbb{P}_{t}^{x,y}(e^{-\left\langle \widehat{l}%
,\chi\right\rangle })=\exp(t(L-M_{_{\frac{\chi}{\lambda}}}))_{x,y}$.

Hence $\mu^{x,y}(e^{-\left\langle \widehat{\gamma},\chi\right\rangle
})=\frac{1}{\lambda_{y}}((I-P+M_{\chi/m})^{-1})_{x,y}=[G_{\chi}]^{x,y}$.

Also $\mathbb{E}^{x}(e^{-\left\langle \widehat{\gamma},\chi\right\rangle
})=\sum_{y}[G_{\chi}]^{x,y}\kappa_{y}$.

\bigskip

In the case of a lattice, one can consider a Poisson process of loops with
intensity $\mu_{00}^{\#}$

\section{Associated Gaussian field}

By a well known calculation, if $X$ is finite, for any $\chi\in\mathbb{R}%
_{+}^{X}$,
\[
\frac{\det(M_{\lambda}-C)}{(2\pi)^{\left|  X\right|  }}\int(e^{-\frac{1}%
{2}<z\overline{z},\chi>}e^{-\frac{1}{2}e(z)}\Pi_{u\in X}\frac{i}{2}%
dz_{u}\wedge d\overline{z}_{u}=\frac{\det(G_{\chi})}{\det(G)}%
\]
and%
\[
\frac{\det(M_{\lambda}+M_{\chi}-C)}{(2\pi)^{\left|  X\right|  }}\int
z^{x}\overline{z}^{y}(e^{-\frac{1}{2}<z\overline{z},\chi>}e^{-\frac{1}{2}%
e(z)}\Pi_{u\in X}\frac{i}{2}dz^{u}\wedge d\overline{z}^{u}=(G_{\chi})^{x,y}%
\]

This can be easily reformulated by introducing the complex Gaussian field
$\phi$ defined by the covariance $\mathbb{E}_{\phi}\mathbb{(}\phi^{x}%
\overline{\phi}^{y})=2G^{x,y}$ (this reformulation cannot be dispensed with
when $X$ becomes infinite)

So we have $\mathbb{E(}(e^{-\frac{1}{2}<\phi\overline{\phi},\chi>}%
)=\det(I+GM_{_{\chi}})^{-1}=\det(G_{\chi}G^{-1})$ and

$\mathbb{E(}(\phi^{x}\overline{\phi}^{y}e^{-\frac{1}{2}<\phi\overline{\phi
},\chi>})=(G_{\chi})^{x,y}\det(G_{\chi}G^{-1})$ Then the following holds:

\begin{theorem}
a) The fields $\widehat{\mathcal{L}_{1}}$ and $\frac{1}{2}\phi\overline{\phi}$
have the same distribution.\label{xy}

\ b) $\mathbb{E}_{\phi}\mathbb{(}(\phi^{x}\overline{\phi}^{y}F(\phi
\overline{\phi}))=\int\mathbb{E}(F(\widehat{\mathcal{L}_{1}}+\widehat{\gamma
}))\mu^{x,y}(d\gamma)$ for any functional $F$ of a non negative field.
\end{theorem}

This is a version of Dynkin's isomorphism (Cf \cite{Dy}). It can be extended
to non symmetric generators (Cf \cite{NS}).

Note it implies immediately that the process $\phi\overline{\phi}$ is
infinitely divisible. See \cite{EK} and its references for a converse and
earlier proofs of this last fact.

In fact an analogous result can be given when $\alpha$ is any positive half
integer, by using a real scalar or vector valued Gaussian field.

Recall that for any $f\in\mathbb{H}$, the law of $f+\phi$ is absolutely
continuous with respect to the law of $\phi$, with density $\exp
(<-Lf,\phi>_{m}-\frac{1}{2}e(f))$

Recall (it was observed by Nelson in the context of the free field) that the
Gaussian field $\phi$ is Markovian: Given any subset $F$ of $X$, denote
$\mathcal{H}_{F}$ the Gaussian space spanned by $\{\phi^{y},y\in F\}$. Then,
for $x\in D=F^{c}$, the projection of $\phi^{x}$ on $\mathcal{H}_{F}$ is
$\sum_{y\in F}[H^{F}]_{y}^{x}\phi^{y}$ .

Moreover, $\phi^{D}=\phi-H^{F}\phi$ is the Gaussian field associated with the
process killed at the exit of $D$.

Note also that if a function $h$ is such that $Lh\leq0$, the loop measure
defined by the $h^{2}m$-symmetric generator $L_{h}=\frac{1}{h}LM_{h}$\ is
associated with the Gaussian field $h\phi$. The killing measure becomes
$\frac{-Lh}{h}\lambda$

Remark finally that the transfer matrix $K$\ is the covariance matrix of the
Gaussian field $d\phi^{x,y}=\phi^{x}-\phi^{y}$ indexed by oriented links.

\section{Energy variation and currents}

The loop measure $\mu$ depends on the energy $e$ which is defined by the free
parameters $C,\kappa$. It will sometimes be denoted $\mu_{e}$. We shall denote
$\mathcal{Z}_{e}$ the determinant $\det(G)=\det(M_{\lambda}-C)^{-1}$. Then
$\mu(p>0)=\log(\mathcal{Z}_{e})+\sum\log(\lambda_{x})$.

Other intrinsic variables of interest on the loop space are associated with
real antisymmetric matrices $\omega_{x,y}$\ indexed by $X^{\Delta}$:
$\omega_{x,y}=-\omega_{y,x}$.. Let us mention a few elementary results.

The operator $[P^{\omega}]_{y}^{x}=P_{y}^{x}\exp(i\omega_{x,y})$ is self
adjoint in $L^{2}(\lambda)$.The associated loop variable writes $\sum
_{j=1}^{p}\omega_{\xi_{j},\xi_{j+1}}$ or $\sum_{x,y}\omega_{x,y}N_{x,y}(l)$.
We will denote it $\int_{l}\omega$. This notation will be used even when
$\omega$ is not antisymmetric. Note it is invariant if $\omega_{x,y}$ is
replaced by $\omega_{x,y}+g(x)-g(y)$ for some $g$. Set $[G^{\omega}%
]^{x,y}=\frac{[(I-P^{\omega})^{-1}]_{y}^{x}}{\lambda_{y}}$ and denote
$\mathcal{Z}_{e,\omega}$ the determinant $\det(G^{\omega})$. By an argument
similar to the one given above for the occupation field, we have:

$\mathbb{P}_{x,x}^{t}(e^{i\int_{l}\omega}-1)=\exp(t(P^{\omega}-I))_{x,x}%
-\exp(t(P-I))_{x,x}$. Integrating in $t$ after expanding, we get from the
definition of $\mu$ :%

\[
\int(e^{i\int_{l}\omega}-1)d\mu(l)=\sum_{k=1}^{\infty}\frac{1}{k}%
[Tr((P^{\omega})^{k})-Tr((P)^{k})]
\]

Hence%
\[
\int(e^{i\int_{l}\omega}-1)d\mu(l)=\log[\det(-L(I-P^{\omega})^{-1}]
\]
and%
\begin{equation}
\mu(\exp(\sum_{l\in\mathcal{L}_{\alpha}}i\int_{l}\omega)-1)=\log
(\det(G^{\omega}G^{-1}))=\log(\frac{\mathcal{Z}_{e,\omega}}{\mathcal{Z}_{e}})
\label{F5}%
\end{equation}

\bigskip

The following result is suggested by an analogy with quantum field theory (Cf
\cite{Gaw}).

\begin{proposition}
i)$\frac{\partial\mu}{\partial\kappa_{x}}=\widehat{l}^{x}\mu$

ii)$\frac{\partial\mu}{\partial\log C_{x,y}}=-T_{x,y}\mu$

with $T_{x,y}(l)=C_{x,y}(\widehat{l}^{x}+\widehat{l}^{y})-N_{x,y}(l)-N_{y,x}(l)$
\end{proposition}

Note that the formula i) would a direct consequence of the Dynkin isomorphism
if we considered only sets defined by the occupation field.

Recall that $\mu=\sum_{x\in X}e^{-\lambda_{x}\tau^{\ast}}\frac{d\tau^{\ast}%
}{\tau^{\ast}}+\sum_{p=2}^{\infty}\sum_{(\xi_{i},i\in\mathbb{Z}/p\mathbb{Z)\in
}X^{p}}\prod_{i\in\mathbb{Z}/p\mathbb{Z}}C_{\xi_{i},\xi_{i+1}}e^{-\lambda
_{\xi_{i}}\tau_{i}^{\ast}}d\tau_{i}^{\ast}$

$C_{x,y}=C_{y,x}=\lambda_{x}P_{y}^{x}$ and $\lambda_{x}=\kappa_{x}+\sum
_{y}C_{x,y}$

The formulas follow by elementary calculation.

Recall that $\mu(\widehat{l}^{x})=G^{x,x}$.and $\mu(N_{x,y})=G^{x,y}C_{x,y}$

So we have $\mu(T_{x,y})=C_{x,y}(G^{x,x}+G^{y,y}-2G^{x,y})$

Then, the above proposition allows to compute all moments of $T$ and
$\widehat{l}$ relative to $\mu_{e}$ (Schwinger functions)

Consider now another energy form $e^{\prime}$ defining an equivalent norm on
$\mathbb{H}$. Then we have the following identity:%
\[
\frac{\partial\mu_{e^{\prime}}}{\partial\mu_{e}}=e^{\sum N_{x,y}%
\log(\frac{C_{x,y}^{\prime}}{C_{x,y}})-\sum(\lambda_{x}^{\prime}-\lambda
_{x})\widehat{l}^{x}}%
\]
The above proposition is the infinitesimal form of this formula. Note that
from the above expression of $\mu$ ($\ref{F2}$),%

\[
\mu_{e}((e^{\sum N_{x,y}\log(\frac{C_{x,y}^{\prime}}{C_{x,y}})-\sum
(\lambda_{x}^{\prime}-\lambda_{x})\widehat{l}^{x}}-1))=\log(\frac{\mathcal{Z}%
_{e^{\prime}}}{\mathcal{Z}_{e}})
\]
(the proof goes by evaluating separately the contribution of trivial loops,
which equals $\sum_{x}\log(\frac{\lambda_{x}}{\lambda_{x}^{\prime}})$).

Note that if $C_{x,y}^{^{\prime}}=h^{x}h^{y}C_{x,y}$ et $\kappa_{x}^{\prime
}=\frac{-Lh}{h}\lambda$ for some positive function $h$ on $E$ such that
$Lh\leq0$, $\frac{\mathcal{Z}_{e^{\prime}}}{\mathcal{Z}_{e}}=\frac{1}%
{\prod(h^{x})^{2}}$.

Note also that$\frac{\mathcal{Z}_{e^{\prime}}}{\mathcal{Z}_{e}}=\mathbb{E(}%
e^{-\frac{1}{2}[e^{\prime}-e](\phi)})$

Equivalently%
\begin{equation}
\mu_{e}(\prod_{(x,y)}[\frac{C_{x,y}^{\prime}}{C_{x,y}}]^{N_{x,y}}\prod
_{x}[\frac{\lambda_{x}}{\lambda_{x}^{\prime}}]^{N_{x}+1}-1)=\mu_{e}%
(\prod_{x,y}[\frac{P_{y}^{\prime x}}{P_{y}^{x}}]^{N_{x,y}}\prod_{x}%
[\frac{\lambda_{x}}{\lambda_{x}^{\prime}}]-1)=\log(\frac{\mathcal{Z}%
_{e^{\prime}}}{\mathcal{Z}_{e}}) \label{F}%
\end{equation}

and therefore%
\[
\mathbb{E}_{\mathbb{E}_{\mathcal{L}_{\alpha}}}(\prod_{(x,y)}[\frac{C_{x,y}%
^{\prime}}{C_{x,y}}]^{N_{x,y}^{(\alpha)}}\prod_{x}[\frac{\lambda_{x}}%
{\lambda_{x}^{\prime}}]^{N_{x}^{(\alpha)}+1})=(\frac{\mathcal{Z}_{e^{\prime}}%
}{\mathcal{Z}_{e}})^{\alpha}%
\]

Note also that $\prod_{(x,y)}[\frac{C_{x,y}^{\prime}}{C_{x,y}}]^{N_{x,y}%
}=\prod_{\{x,y\}}[\frac{C_{x,y}^{\prime}}{C_{x,y}}]^{N_{x,y}+N_{y,x}}$

\bigskip

N.B.: These $\frac{\mathcal{Z}_{e^{\prime}}}{\mathcal{Z}_{e}}$ determine, when
$e^{\prime}$ varies with $\frac{C^{^{\prime}}}{C}\leq1$ and $\frac{\lambda
^{\prime}}{\lambda}=1$, the Laplace transform of the distribution of the
traversal numbers of non oriented links $N_{x,y}+N_{y,x}$, hence the loop
distribution $\mu_{e}$.

More generally%
\begin{equation}
\mu_{e}(e^{-\sum N_{x,y}\log(\frac{C_{x,y}^{^{\prime}}}{C_{x,y}})-\sum
(\lambda_{x}^{^{\prime}}-\lambda_{x})\widehat{l}_{x}+i\int_{l}\omega}%
-1)=\log(\frac{\mathcal{Z}_{e^{^{\prime}},\omega}}{\mathcal{Z}_{e}})
\label{F4}%
\end{equation}
or%
\[
\mu_{e}(\prod_{x,y}[\frac{C_{x,y}^{\prime}}{C_{x,y}}e^{i\omega_{x,y}%
}]^{N_{x,y}}\prod_{x}[\frac{\lambda_{x}}{\lambda_{x}^{\prime}}]^{N_{x}%
+1}-1)=\log(\frac{\mathcal{Z}_{e^{\prime},\omega}}{\mathcal{Z}_{e}})
\]
Note also that this last formula applies to the calculation of loop indices if
we have for exemple a simple random walk on an oriented two dimensional
lattice. In such cases, $\omega_{z^{\prime}}$ can be chosen such that
$\int_{l}\omega_{z^{\prime}}$\ is the winding number of the loop around a
given point $z^{\prime}$\ of the dual lattice$\footnote{The construction of
$\omega$ can be done as follows: Let $P^{\prime}$ be the uniform Markov
transition probability on neighbouring points of the dual lattice and let $h$
be a function such that $P^{\prime}h=h$ except in $z^{\prime}$. Then if the
link $xy$ in $X$ intersects $x^{\prime}y^{\prime}$ in $X^{\prime}$, with
$\det(x-y,x^{\prime}-y^{\prime})>0$, set $\omega_{x,y}=h(y^{\prime
})-h(x^{\prime})$}$ $X^\prime$. Then $e^i\pi\sum_l\in\mathcal{L}_\alpha
\int_l\omega_z^\prime$ is a spin system of interest.

We then get for exemple that%
\[
\mu(\int_{l}\omega)\neq0)=-\frac{1}{2\pi}\int_{0}^{2\pi}\log(\det(G^{2\pi
u\omega}G^{-1}))du
\]
and hence%
\[
\mathbb{P(}\sum_{l\in\mathcal{L}_{\alpha}}|\int_{l}\omega_{z}^{\prime
}|)=0)=e^{\frac{\alpha}{2\pi}\int_{0}^{2\pi}\log(\det(G^{2\pi u\omega}%
G^{-1}))du}%
\]

Conditional distributions of the occupation field with respect to values of
the winding number can also be obtained.

We can apply the formula \ref{F}\ to calculations concerning the links visited
by the loops (similar to those done in section \ref{occup}\ for sites).

For exemple, $R$ is a set of links, denote $e^{\left]  R\right[  }$ the energy
form defined \ from $e$ by setting all conductances in $R$ to zero and
increasing $\kappa$ in such a way that $\lambda$ is unchanged..

Then $\mu_{e}(\sum_{(x,y)\in R}N_{x,y}+N_{y,x}>0)=-\log(\frac{\det(G^{\left]
R\right[  })}{\det(G)})$ and therefore, the probability no loop in
$\mathcal{L}_{\alpha}$ visits $R$ equals $\frac{\det(G^{\left]  R\right[  }%
)}{\det(G)}=(\frac{\mathcal{Z}_{e^{\left]  R\right[  }}}{\mathcal{Z}_{e}%
})^{\alpha}$.

\section{Self-avoiding paths and spanning trees.}

Recall that link $f$ is a pair of points $(f^{+},f^{-})$ such that
$C_{f}=C_{f^{+},f^{-}}\neq0$. Define $-f=(f^{-},f^{+})$.

Let $\mu_{x,y}^{\neq}$ be the measure induced by $C$ on discrete self-avoiding
paths.between $x$\ and $y$: $\mu_{\neq}^{x,y}(x,x_{2},...,x_{n-1}%
,y)=C_{x,x_{2}}C_{x_{1},x_{3}}...C_{x_{n-1},y}$.

Another way to defined a measure on discrete self avoiding paths from $x$ to
$y$ is loop erasure (see for exemple \cite{Law}). One checks easily the following:

\begin{proposition}
the image of $\mu^{x,y}$ by the loop erasure map $\gamma\rightarrow\gamma
^{BE}$\ is $\mu_{BE}^{x,y}$ defined on self avoiding paths by $\mu_{BE}%
^{x,y}(\eta)=\mu_{\neq}^{x,y}(\eta)\frac{\det(G)}{\det(G^{\{\eta\}^{c}})}%
=\mu_{\neq}^{x,y}(\eta)\det(G_{|\{\eta\}\times\{\eta\}})$ (Here $\{\eta\}$
denotes the set of points in the path $\eta$)
\end{proposition}

Proof: If $\eta=(x_{1}=x,x_{2},...x_{n}=y)$,and $\eta_{m}=(x,...x_{m})$,
$\mu^{x,y}(\gamma^{BE}=\eta)=V_{x}^{x}P_{x_{2}}^{x}[V^{\{x\}^{c}}]_{x_{2}%
}^{x_{2}}...[V^{\{\eta_{n-1}\}^{c}}]_{x_{n-1}}^{x_{n-1}}P_{y}^{x_{n-1}%
}[V^{\{\eta\}^{c}}]_{y}^{y}\lambda_{y}^{-1}=\mu_{\neq}^{x,y}(\eta
)\frac{\det(G)}{\det(G^{\{\eta\}^{c}})}$ as $[V^{\{\eta_{m-1}\}^{c}}]_{x_{m}%
}^{x_{m}}=\frac{\det([(I-P]|_{\{\eta_{m}\}^{c}\times\{\eta_{m}\}^{c}})}%
{\det([(I-P]|_{\{\eta_{m-1}\}^{c}\times\{\eta_{m-1}\}^{c}})}=\frac{\det
(V^{\{\eta_{m-1}\}^{c}})}{\det(V^{\{\eta_{m}\}^{c}})}=\frac{\det
(G^{\{\eta_{m-1}\}^{c}})}{\det(G^{\{\eta_{m}\}^{c}})}\lambda^{x_{m}\text{.}}$
for all $m\leq n-1$.

Also:$\int$ $e^{-<\widehat{\gamma},\chi>}1_{\{\gamma^{BE}=\eta\}}\mu
^{x,y}(d\gamma)=\frac{\det(G_{\chi})}{\det(G_{\chi}^{\{\eta\}^{c}}%
)}e^{-<\widehat{\eta},\chi>}\mu_{\neq}^{x,y}(\eta)$

$=\det(G_{\chi})_{|\{\eta\}\times\{\eta\}}e^{-<\widehat{\eta},\chi>}\mu_{\neq
}^{x,y}(\eta)=$ $\frac{\det(G_{\chi})_{|\{\eta\}\times\{\eta\}}}%
{\det(G_{|\{\eta\}\times\{\eta\}})}e^{-<\widehat{\eta},\chi>}\mu_{BE}%
^{x,y}(\eta)$ for any self-avoiding path $\eta$.

Therefore, under $\mu^{x,y}$, the conditional distribution of $\widehat
{\gamma}-\widehat{\eta}$ given $\gamma^{BE}=\eta$ is the distribution of
$\widehat{\mathcal{L}_{1}}-\widehat{\mathcal{L}_{1}^{\{\eta\}^{c}}\text{ }}$
i.e. the occupation field of the loops of $\mathcal{L}_{1}$\ which intersect
$\eta$.

More generally, it can be shown that

\begin{proposition}
\label{be}the conditional distribution of the set $\mathcal{L}_{\gamma}%
$\ of\ loops of $\gamma$ given $\gamma^{BE}=\eta$ is the distribution of
$\mathcal{L}_{1}/\mathcal{L}_{1}^{\{\eta\}^{c}}$ i.e. the loops of
$\mathcal{L}_{1}$\ which intersect $\eta$.
\end{proposition}

Proof: First an elementary calculation shows that

$\mu_{e^{\prime}}^{x,y}(\gamma^{BE}=\eta)=\frac{C_{x,x_{2}}^{\prime}%
C_{x_{1},x_{3}}^{\prime}...C_{x_{n-1},y}^{\prime}}{C_{x,x_{2}}C_{x_{1},x_{3}%
}...C_{x_{n-1},y}}\mu_{e}^{x,y}(\prod_{u\neq v}[\frac{C_{u,v}^{\prime}%
}{C_{u,v}}]^{N_{x,y}(\mathcal{L}_{\gamma})+N_{y,x}(\mathcal{L}_{\gamma})}%
\prod_{u}[\frac{\lambda_{u}}{\lambda_{u}^{\prime}}]^{N_{u}(\mathcal{L}%
_{\gamma})}1_{\{\gamma^{BE}=\eta\}})$

Therefore, by the previous proposition,

$\mu_{e}^{x,y}(\prod_{u\neq v}[\frac{C_{u,v}^{\prime}}{C_{u,v}}]^{N_{x,y}%
(\mathcal{L})+N_{y,x}(\mathcal{L})}\prod_{u}[\frac{\lambda_{u}}{\lambda
_{u}^{\prime}}]^{N_{u}(\mathcal{L})}|\gamma^{BE}=\eta)=\frac{\mathcal{Z}%
_{e}\mathcal{Z}_{[e^{\prime}]^{\{\eta\}^{c}}}}{\mathcal{Z}_{e^{\{\eta\}^{c}}%
}\mathcal{Z}_{e^{\prime}}}$.

Moreover, by \ref{F} and the properties of the Poisson processes,

$\mathbb{E}(\prod_{u\neq v}[\frac{C_{u,v}^{\prime}}{C_{u,v}}]^{N_{x,y}%
(\mathcal{L}_{1}/\mathcal{L}_{1}^{\{\eta\}^{c}})+N_{y,x}(\mathcal{L}%
_{1}/\mathcal{L}_{1}^{\{\eta\}^{c}})}\prod_{u}[\frac{\lambda_{u}}{\lambda
_{u}^{\prime}}]^{N_{u}(\mathcal{L}_{1}/\mathcal{L}_{1}^{\{\eta\}^{c}})}=$
$\frac{\mathcal{Z}_{e}\mathcal{Z}_{[e^{\prime}]^{\{\eta\}^{c}}}}%
{\mathcal{Z}_{e^{\{\eta\}^{c}}}\mathcal{Z}_{e^{\prime}}}$

It follows that the distributions of the $N_{x,y}+N_{y,x}$'s are identical for
the set of erased loops and $\mathcal{L}_{1}/\mathcal{L}_{1}^{\{\eta\}^{c}}$.
Moreover, remark \ref{unif} allows to conclude, since the same conditional
equidistribution property holds for the configurations of erased loops.

Similarly one can define the image of $\mathbb{P}^{x}$ by $BE$ which is given by

$\mathbb{P}_{BE}^{x}(\eta)=C_{x_{1},x_{2}}...C_{x_{n-1},x_{n}}\kappa_{x_{n}%
}\det(G_{|\{\eta\}\times\{\eta\}})$, for $\eta=(x_{1},...,x_{n})$, and get the
same results.

\bigskip

Wilson's algorithm (see \cite{Lyo2}) iterates this construction, starting with
$x^{\prime}s$ in arbitrary order. Each step of the algorithm reproduces the
first step except it stops when it hits the already constructed tree of self
avoiding paths. It provides a construction of the probability measure
$\mathbb{P}_{ST}^{e}$\ on the set $ST_{X,\Delta}$\ of spanning trees of $X$
rooted at the cemetery point $\Delta$ defined by the energy $e$. The weight
attached to each oriented link $\xi=(x,y)$ of $X\times X$ is the conductance
and the weight attached to the link $(x,\Delta)$ is $\kappa_{x}$. As the
determinants simplify, the probability of a tree $\Upsilon$ is given by the
simple formula%

\[
\mathbb{P}_{ST}^{e}(\Upsilon)=\mathcal{Z}_{e}\prod_{\xi\in\Upsilon}C_{\xi}%
\]

\begin{proposition}
\label{wa}The random set of discrete loops $\mathcal{L}_{W}$\ constructed in
this algorithm is independent of the random spanning tree, and independent of
the ordering. It has the same Poisson distribution as the non trivial discrete
loops of $\mathcal{L}_{1}$.
\end{proposition}

\bigskip It follows easily from proposition \ref{be} .

Together with the spanning tree these discrete loops define an interesting
random graph$.$

\bigskip

First note that, since we get a probability
\[
\mathcal{Z}_{e}\sum_{\Upsilon\in ST_{X,\Delta}}\prod_{(x,y)\in\Upsilon}%
C_{x,y}\prod_{x,(x,\Delta)\in\Upsilon}\kappa_{x}=1
\]
or equivalently
\[
\sum_{\Upsilon\in ST_{X,\Delta}}\prod_{(x,y)\in\Upsilon}P_{y}^{x}%
\prod_{x,(x,\Delta)\in\Upsilon}P_{\Delta}^{x}=\frac{1}{\prod_{x\in X}%
\lambda_{x}\mathcal{Z}_{e}}%
\]

so that
\[
\mathbb{P}_{ST}^{e}(\Upsilon)=\mathcal{Z}_{e}\prod_{x\in X}\lambda_{x}%
\prod_{(x,y)\in\Upsilon}P_{y}^{x}\prod_{x,(x,\Delta)\in\Upsilon}P_{\Delta}^{x}%
\]

Then, it comes that, for any $e^{\prime}$,%

\[
\mathbb{E}_{ST}^{e}(\prod_{(x,y)\in\Upsilon}\frac{P_{y}^{\prime x}}{P_{y}^{x}%
}\prod_{x,(x,\Delta)\in\Upsilon}\frac{P_{\Delta}^{\prime x}}{P_{\Delta}^{x}%
})=\frac{\prod_{x\in X}\lambda_{x}}{\prod_{x\in X}\lambda_{x}^{\prime}%
}\frac{\mathcal{Z}_{e}}{\mathcal{Z}_{e^{\prime}}}%
\]
and%

\begin{equation}
\mathbb{E}_{ST}^{e}(\prod_{(x,y)\in\Upsilon}\frac{C_{x,y}^{\prime}}{C_{x,y}%
}\prod_{x,(x,\Delta)\in\Upsilon}\frac{\kappa_{x}^{\prime}}{\kappa_{x}%
})=\frac{\mathcal{Z}_{e}}{\mathcal{Z}_{e^{\prime}}} \label{Q}%
\end{equation}

We also have $\mathbb{P}_{ST}^{e}((x,y)\in\Upsilon)=\mathbb{P}_{x}^{BE}%
(\eta_{1}=y)=V_{x}^{x}P_{y}^{x}\mathbb{P}^{y}(T_{x}=\infty)=C_{x,y}%
G^{x,x}(1-\frac{G^{x,y}}{G^{x,x}})$

From the results exposed in \cite{Lyo} and \cite{Lyo2}, or directly from the
above, we recover Kirchhoff's theorem:

$\mathbb{P}_{ST}^{e}(\pm(x,y)\in\Upsilon)=C_{x,y}[G^{x,x}(1-\frac{G^{x,y}%
}{G^{x,x}})+G^{y,y}(1-\frac{G_{x,y}}{G_{y,y}})]=C_{x,y}(G^{x,x}+G^{y,y}%
-2G^{x,y})=C_{x,y}K^{x,y),(x,y)}$ and more generally Pemantle's transfer
current theorem:%

\[
\mathbb{P}_{ST}^{e}(\pm\xi_{1},...\pm\xi_{k}\in\Upsilon)=(\prod_{1}^{k}%
C_{\xi_{i}})\det(K^{\xi_{i},\xi_{j}}\;1\leq i,j\leq k)
\]

Note this determinant does not depend on the orientation of the links.

\bigskip

Proof: We use recurrence on $k$. Let $M$ denote the smallest subset of
$X^{\Delta}$ containing the links $\pm\xi_{1},...\pm\xi_{k}$ and denote
$E-M$\ by $D$ Let $V$ be the subspace of $\mathbb{A}$ spanned by all
$K^{(x,y)}$ with $x$ and $y$ in $M_{k}$. Note that the orthogonal of $V$ in
$\mathbb{A}$ is spanned by $dG^{D}(\delta_{x})$ and that for any $\eta=(u,v)$
the projection of $K^{\eta}$ on $V^{\bot}$ is $dG^{D}(\delta_{v}-\delta_{u})$,
and $\left\langle dG^{D}(\delta_{v}-\delta_{u}),dG^{D}(\delta_{v}-\delta
_{u})\right\rangle =[G^{D}]^{u,u}+[G^{D}]^{v,v}-2[G^{D}]^{u,v}$

Moreover $\det(K^{\xi_{i},\xi_{j}}\;1\leq i,j\leq k)=\left\|  K^{\xi_{1}%
}\wedge...\wedge K^{\xi_{k}}\right\|  _{\wedge^{k}\mathbb{A}}^{2}$. Therefore,
if $\xi_{k+1}=\eta$ $\det(K^{\xi_{i},\xi_{j}}\;1\leq i,j\leq k+1)=\det
(K^{\xi_{i},\xi_{j}}\;1\leq i,j\leq k)([G^{D}]^{u,u}+[G^{D}]^{v,v}%
-2[G^{D}]^{u,v})$

But the argument given for $k=1$ shows also that $\mathbb{P}_{ST}^{e}(\pm
\eta\in\Upsilon|\pm\xi_{1},...\pm\xi_{k}\in\Upsilon)=C_{u,v}([G^{D}%
]^{u,u}+[G^{D}]^{v,v}-2[G^{D}]^{u,v})$ so we can conclude.

\bigskip

Therefore, given any function $g$ on non oriented links, $\mathbb{E}_{ST}%
^{e}(e^{-\sum_{\xi\in\Upsilon}g(\xi)})=\mathbb{E}_{ST}^{e}(\prod_{\xi
}(1+(e^{-g(\xi)}-1)1_{\xi\in\Upsilon})=\sum Tr((M_{C(e^{-g}-1)}K)^{\wedge k})$
and we have
\[
\mathbb{E}_{ST}^{e}(e^{-\sum_{\xi\in\Upsilon}g(\xi)})=\det(I+KM_{C(e^{-g}%
-1)})=\det(I-M_{\sqrt{C(1-e^{-g})}}KM_{\sqrt{C(1-e^{-g})}})
\]
\textsl{ }This is an exemple of the Fermi point processes discussed in
\cite{ShiTak}\textsl{.}

\textsl{ }But, by (\ref{Q}) and (\ref{F}), it comes that%

\[
\log(\mathbb{E}_{ST}^{e}(\prod_{(x,y)\in\Upsilon}\frac{C_{x,y}^{\prime}%
}{C_{x,y}}\prod_{x,(x,\Delta)\in\Upsilon}\frac{\kappa_{x}^{\prime}}{\kappa
_{x}})=-\mathbb{\mu}_{e}(\prod_{x,y}[\frac{C_{x,y}^{\prime}}{C_{x,y}%
}]^{N_{x,y}}\prod_{x}[\frac{\lambda_{x}}{\lambda_{x}^{\prime}}]^{N_{x}%
+1}-1)=\log(\frac{\mathcal{Z}_{e}}{\mathcal{Z}_{e^{\prime}}})
\]

The first identity could also be derived from proposition.\ref{wa}$.$

As $\frac{\lambda_{x}}{\lambda_{x}^{\prime}}=\frac{1}{\sum P_{y}%
^{x}\frac{C_{x,y}^{\prime}}{C_{x,y}}+P_{\Delta}^{x}\frac{\kappa_{x}^{\prime}%
}{\kappa_{x}}}$, (with the convention $\frac{P_{\Delta}^{x}}{\kappa_{x}%
}=\lambda_{x}^{-1}$ if $\kappa_{x}=0$)\ we obtain

For any function $g$ on non oriented link of $X_{\Delta},$ non negative on
links of $X$
\[
\mathbb{E}_{\mathcal{L}_{\alpha}}(e^{-\sum_{x,y}g(\{x,y\})N_{x,y}}\prod
_{x}[\sum_{z\in X^{\Delta}}P_{z}^{x}e^{-g\{x,z\}}]^{-N_{x}-1})=\det
(I-KM_{C(1-e^{-g})})^{-\alpha}%
\]

We can check that this formula allows to recover the identity $\mathbb{E}%
_{\mathcal{L}_{\alpha}}(\frac{N_{x}+N_{y}+2}{\lambda_{x}}-N_{x,y}%
-N_{y,x})=\alpha C_{x,y}(G^{x,x}+G^{y,y}-2G^{x,y})$ It also gives back
prpositin \ref{det1}\ for $g=0$ on $X\times X$.

If $\kappa$ is positive everywhere, we can adjust $g(\{x,\Delta\})$ to make
$\sum_{z\in X^{\Delta}}P_{z}^{x}e^{-g\{x,z\}}=1$ This means we have to choose
$\kappa_{x}(1-e^{-g(\{x,\Delta\})})=\kappa_{x}(1-\frac{\lambda_{x}}{\kappa
_{x}}(1-\sum_{z\in X}P_{z}^{x}e^{-g(\{x,z\})}))$

$=(\kappa_{x}-\lambda_{x}+$.$\sum_{z\in X}C_{xz}e^{-g(\{x,z\})})=\sum_{z\in
X}C_{xz}(e^{-g(\{x,z\})}-1)$

We check also also that by \ref{F}

$\mathbb{E}_{\mathcal{L}_{\alpha}}(e^{-\sum_{x,y}g(\{x,y\})N_{x,y}%
})=(\frac{\mathcal{Z}_{e^{\prime}}}{\mathcal{Z}_{e}})^{\alpha}=(\frac{\det
(M_{\lambda}-Ce^{-g})}{\det(M_{\lambda}-C)})^{-a}$ Finally, the restriction on
$\kappa$ can be removed by taking a limit and we obtain:

\begin{proposition}
For any function $g$ on non oriented link of $X_{\Delta},$ non negative on
links of $X$, set $Tg(\xi)=C_{\xi}(1-e^{-g(\xi)})$ if $\xi$ is a link of $X$
and $Tg(\{x,\Delta\})=\sum_{z\in X}C_{xz}(e^{-g(\{x,z\})}-1)$ for all $x$.
Then%
\[
\det(I+C(I-[e^{-g}]))^{-\alpha}=\mathbb{E}_{\mathcal{L}_{\alpha}}%
(e^{-\sum_{x,y}g(\{x,y\})N_{x,y}})=\det(I-KM(g))^{-\alpha}%
\]
\end{proposition}

\bigskip

We see that the Poisson measure on loops $\mathcal{L}_{\alpha}$ induces a
point process $N$\ on the space of non oriented links defined by the pair
($\alpha,K)$ which reminds the point processes discussed in \cite{ShiTak}.
Note however a difference of sign in the right hand side determinant, which is
not a Laplace transform for positive $\alpha$.

\section{Fock spaces and Wick product}

Recall that the Gaussian space $\mathcal{H}$ spanned by $\{\phi^{x},x\in X\}$
is isomorphic to $\mathbb{H}$ by the linear map mapping $\operatorname{Re}%
(\phi_{x})$ on $G_{x,\cdot}$ which extends into an isomorphism between \ the
space of square integrable functionals of the Gaussian fields and the
symmetric Fock space obtained as the closure of the sum of all symmetric
tensor powers of $\mathbb{H}$ (Bose second quantization). We have seen that
$L^{2}$ functionals of $\widehat{\mathcal{L}_{1}}$ can be represented in this
symmetric Fock space.

In order to prepare the extension of these isomorphisms to a more interesting
framework (including especially the planar Brownian motion considered in
\cite{LW}) we shall introduce the renormalized (or Wick) powers of
$\phi\overline{\phi}$.

The Laguerre polynomials $L_{n}^{0}(x)$ are defined by their generating
function $\frac{e^{-\frac{xt}{1-t}}}{1-t}=\sum t^{n}L_{n}^{0}(x)$

Then one defines the polynomial $P_{n}(\cdot)=(-1)^{n}n!L_{n}^{0}(\cdot)$

Setting $\sigma_{x}=G_{x,x}$, it comes that $\sigma^{n}P_{n}(\frac{\phi
\overline{\phi}}{\sigma})$ is the inverse image of a $2n$-th tensor\ in the
Fock space denoted $:(\phi\overline{\phi})^{n}:$.Note that $:(\phi
\overline{\phi}):=(\phi\overline{\phi})-2\sigma$ These variables are
orthogonal in $L^{2}$. Set $\widetilde{l}^{x}=\widehat{l}^{x}-\sigma_{x}$ be
the centered occupation field. Note that an equivalent formulation of
proposition \ref{xy} is that the fields $\frac{1}{2}$ $:\phi\overline{\phi}:$
and $\widetilde{\mathcal{L}}_{1}$\ have the same law.

Let us now consider the relation of higher Wick powers with self intersection
local times.

\section{Decompositions}

If $D\subset X$ and we set $F=D^{c}$, the orthogonal decomposition of the
Dirichlet norm $e(f)$\ into $e^{D}(f-H^{F}f)+e(H^{F}f)$ (cf \cite{LJ0} and
references) leads to the decomposition of the Gaussian field mentionned above
and also to a decomposition of the Markov chain into the Markov chain killed
at the exit of $D$ and the trace of the Markov chain on $F$.

\begin{proposition}
The trace of the Markov chain on $F$ is defined by the Dirichlet norm
$e^{\{F\}}(f)=e(H^{F}f)$ , for which
\[
C_{x,y}^{\{F\}}=C_{x,y}+\sum_{a,b\in D}C_{x,a}C_{b,y}[G^{D}]^{a,b}%
\]
\[
\lambda_{x}^{\{F\}}=\lambda_{x}-\sum_{a,b\in D}C_{x,a}C_{b,x}[G^{D}]^{a,b}%
\]
and%
\[
\mathcal{Z}_{e}=\mathcal{Z}_{e^{D}}\mathcal{Z}_{e^{\{F\}}}%
\]
\end{proposition}

Proof: The first assertion is well known. For the second, note first that for
any $y\in F$, $[H^{F}]_{y}^{x}=1_{x=y}+1_{D}(x)\sum_{b\in D}[G^{D}%
]^{x,b}C_{b,y}$. Moreover, $e(H^{F}f)=\left\langle f,H^{F}f\right\rangle _{e}$
and therefore

$\lambda_{x}^{\{F\}}=e^{\{F\}}(1_{\{x\}})=e(1_{\{x\}},H^{F}1_{\{x\}}%
)=\lambda_{x}-\sum_{a\in D}C_{x,a}[H^{F}]_{x}^{a}=\lambda_{x}-\sum_{a,b\in
D}C_{x,a}C_{b,x}[G^{D}]^{a,b}$.

Then for distinct $x$ and $y$ in $F$,

$C_{x,y}^{\{F\}}=-\left\langle 1_{\{x\}},1_{\{y\}}\right\rangle _{e^{\{F\}}%
}=-\left\langle 1_{\{x\}},H^{F}1_{\{y\}}\right\rangle _{e}=C_{x,y}+\sum
_{a}C_{x,a}[H^{F}]_{y}^{a}=C_{x,y}+$ $\sum_{a,b\in D}C_{x,a}C_{b,y}%
[G^{D}]^{a,b}$.

Finally, note also that $G^{\{F\}}$ is the restriction of $G$ to $F$. as for
all $x,y\in F$, $\left\langle G^{\cdot,y},1_{\{x\}}\right\rangle _{e^{\{F\}}%
}=\left\langle G^{\cdot,y},[H^{F}]_{x}^{\cdot}\right\rangle _{e}=1_{\{x=y\}}$.
Hence the determinant decomposition already used in yields the final formula.

The cases where $F$ has one point was already treated in section 3-2.

The transition matrix $[P^{\{F\}}]_{y}^{x}$ can also be computed directly and equals

$P_{y}^{x}+$ $\sum_{a,b\in D}P_{a}^{x}P_{y}^{b}V^{D\cup\{x\}}]_{b}^{a}%
=P_{y}^{x}+$ $\sum_{a,b\in D}P_{a}^{x}C_{b,y}[G^{D\cup\{x\}}]^{a,b}$. The
calculation of$\frac{C_{x,y}^{\{F\}}}{\lambda_{x}^{\{F\}}}$ yields a
decomposition in two parts according whether the jump to $y$ occurs from $x$
or from $D$.

\bigskip

If we set $e_{\chi}=e+\left\|  {}\right\|  _{L^{2}(\chi)}$ and denote
$[e_{\chi}]^{\{F\}}$ by $e^{\{F,\chi\}}$ we have
\[
C_{x,y}^{\{F,\chi\}}=C_{x,y}+\sum_{a,b}C_{x,a}C_{b,y}[G_{\chi}^{D}]^{a,b}%
\]

and
\[
\lambda_{x}^{\{F,\chi\}}=\lambda_{x}-\sum_{a,b}C_{x,a}C_{b,x}[G_{\chi}%
^{D}]^{a,b})
\]

\bigskip More generally, if $e^{\#}$ is such that $C^{\#}=C$\ on $F\times F$,
and $\lambda=\lambda^{\#}$ on $F$ we have:
\[
C_{x,y}^{\#\{F\}}=C_{x,y}+\sum_{a,b}C_{x,a}^{\prime}C_{b,y}^{\prime}[G_{\chi
}^{\prime D}]^{a,b}%
\]

and
\[
\lambda_{x}^{\#\{F\}}=\lambda_{x}-\sum_{a,b}C_{x,a}^{\prime}C_{b,x}^{\prime
}[G_{\chi}^{\prime D}]^{a,b})
\]

\bigskip

A loop in $X$ which hits $F$\ can be decomposed into a loop $l^{\{F\}}$\ in
$F$ and its excursions in $D$ which may come back to their starting point.

Set $\nu_{x,y}^{D}=C_{x,y}\delta_{\emptyset}+\sum_{a,b\in D}C_{x,a}C_{b,y}%
\mu_{D}^{a,b}$ and $\nu_{D}^{x}=\lambda_{x}^{-1}\delta_{\emptyset}+\sum
_{n=1}^{\infty}\lambda_{x}^{-n}(\sum_{a,b\in D}C_{x,a}C_{b,x}\mu_{D}%
^{a,b})^{\otimes n}$. Here $\mu_{D}^{a,b}$\ denotes the bridge measure (with
mass $[G^{D}]^{a,b}$\ associated with $e^{D}$. Note that $\nu_{x,y}%
^{D}(1)=C_{x,y}^{\{F\}}$ and $\nu_{x}^{D}(1)=\frac{1}{\lambda_{x}^{\{F\}}}$.
We get a decomposition of $\mu$ into its restriction $\mu^{D}$\ to loops in
$D$ (associated to the process killed at the exit of $D$), a loop measure
$\mu^{\{F\}}$\ defined on loops of $F$ by the trace of the Markov chain on
$F$, measures $\nu_{x,y}^{D}$\ on excursions in $D$ indexed by pairs of points
in $F$ and measures $\nu_{x}^{D}\ $on finite sequences of excursions in
$D$\ indexed by points of $F$.

Conversely, a loop $l^{\{F\}}$ of points $\xi_{i}$\ in $F$ (possibly reduced
to a point), a family of excursions $\gamma_{\xi_{i},\xi_{i+1}}$\ attached to
the jumps of $l^{\{F\}}$ and systems of i.i.d. excursions $,\gamma_{\xi_{i}}$
attached to the points of $l^{\{F\}}$ defines a loop $\Lambda(l^{\{F\}}%
,(\gamma_{\xi_{i},\xi_{i+1}}),(\gamma_{\xi_{i}}))$. Note excursions can be
empty. Then $\mu-\mu^{D}$ is the image measure by $\Lambda$ of $\mu
^{\{F\}}(dl^{\{F\}})\prod(\widetilde{\nu}_{\xi_{i},\xi_{i+1}}^{D}%
)(d\gamma_{\xi_{i},\xi_{i+1}})\prod\widetilde{\nu}_{\xi_{i}}^{D}(d\gamma
_{\xi_{i}})$, denoting $\widetilde{\nu}$ the normalised measures $\frac{\nu
}{\nu(1)}$.

The Poisson process $\mathcal{L}_{\alpha}^{\{F\}}=\{l^{\{F\}},l\in
\mathcal{L}_{\alpha}\}$ has intensity $\mu^{\{F\}}$ and is independent of
$\mathcal{L}_{\alpha}^{D}$.

In particular, if $\chi$ is a measure carried by $D$, we have:%

\[
E(e^{-\left\langle \widehat{\mathcal{L}_{\alpha}},\chi\right\rangle
}|\mathcal{L}_{\alpha}^{\{F\}})=[\frac{\mathcal{Z}_{e_{\chi}^{D}}}%
{\mathcal{Z}_{e^{D}}}]^{\alpha}\prod_{l^{\{F\}}\in\mathcal{L}_{\alpha}%
^{\{F\}}}(\prod_{x,y\in F}[\frac{C_{x,y}^{\{F,\chi\}}}{C_{x,y}^{\{F\}}%
}]^{N_{x,y}}\prod_{x\in F}[\frac{\lambda_{x}^{\{F\}}}{\lambda_{x}^{\{F,\chi
\}}}]^{N_{x}+1})
\]

More generally

\begin{proposition}
if $C^{\#}=C$\ on $F\times F$, and $\lambda=\lambda^{\#}$ on $F$

$\mathbb{E}_{\mathcal{L}_{\alpha}}(\prod_{x,y\notin F\times F}[\frac{C_{x,y}%
^{\#}}{C_{x,y}}]^{N_{x,y}}e^{-\sum_{x\in D}\widehat{l_{x}}(\lambda_{x}%
^{\#}-\lambda_{x})}|\mathcal{L}_{\alpha}^{\{F\}})=$

$[\frac{\mathcal{Z}_{e^{\#D}}}{\mathcal{Z}_{e^{D}}}]^{\alpha}\prod
_{l^{\{F\}}\in\mathcal{L}_{\alpha}^{\{F\}}}(\prod_{x,y\in F}[\frac{C_{x,y}%
^{\#\{F\}}}{C_{x,y}^{\{F\}}}]^{N_{x,y}}\prod_{x\in F}[\frac{\lambda
_{x}^{\{F\}}}{\lambda_{x}^{\prime\{F\}}}]^{N_{x}+1})$
\end{proposition}

The proof can be done by decomposing all $e^{\prime}$ into $e^{\#}+(e^{\prime
}-e^{\#})$, with $(e^{\prime}-e^{\#})$ carried by$F$

\bigskip

These decomposition formulas extend to include a current $\omega$ provided it
is closed (i.e. vanish on every loop) in $D$. In particular, it allows to
define $\omega^{F}$ such that:%
\[
\mathcal{Z}_{e,\omega}=\mathcal{Z}_{e^{D}}\mathcal{Z}_{e^{\{F\}},\omega^{F}}%
\]

\section{Reflection positivity and Hilbert space}

Let us fix $\alpha$. In view of physical applications, it is appropriate to
assume that $X$ is the union of two parts $X^{\pm}$ exchanged by an involution
$\rho$ under which $e$ is invariant. Each configuration $\mathcal{L}_{\alpha}%
$\ of loops induces a configuration $\Lambda$\ of loops in $X^{0}=X^{+}\cap
X^{-}$. Given a function $F$ on loops configuration in $X^{+}$, it follows
from the previous proposition the following

\begin{corollary}
$E(F(\mathcal{L}_{\alpha}|_{X^{+}})|\Lambda)=E(F\circ\rho(\mathcal{L}_{\alpha
}|_{X^{-}})|\Lambda)$
\end{corollary}

so that \ the reflection positivity (also called physical positivity) property
holds:%
\[
\left\|  F\right\|  =E(F(\overline{F}\circ\rho))\geq0.
\]

The physical space is the quotient space modulo functionals of zero norm. It
identifies with $L^{2}$\ functionals of $\Lambda$. Osterwalder-Schrader-type
construction can be used to produce non commuting field observables. More
precisely, after extending the framework to infinite spaces (see section
below), one can assume for exemple $X$ has a product structure $S\times
\mathbb{Z}$ and that the time translation $\tau$\ and the time reversal $\rho$
leave $e$ invariant. Then $\tau$ induces a self adjoint contraction $T$\ of
the physical space, hence a Hamiltonian $\log(T)$\ and by complex
exponentiation, a unitary dynamic $U$. Non commuting observables are obtained
by conjugation of an observable by the opeators $U^{n}$. This extends the
construction of the relativistic non commuting quantum free field observables
out of the Euclidean Gaussian field.

\section{The case of general Markov processes}

We now explain briefly how some of the above results will be extended to a
symmetric Markov process on an infinite space $X$. The construction of the
loop measure as well as a lot of computations can be performed quite
generally, using Dirichlet space theory \ Let us consider more closely the
occupation field $\widehat{l}$. \ The extension is rather straightforward when
points are not polar. We can start with a Dirichlet space of continuous
functions and a measure $m$ such that there is a mass gap. Let $P_{t}$ the
associated Feller semigroup. Then the Green function is well defined as the
mutual energy of the Dirac measures $\delta_{x}$ and $\delta_{y}$ which have
finite energy. It is the covariance function of a Gaussian Markov field
$\phi^{x}$, which will be associated to the field $\widehat{l}$ of local
times.of the Poisson process of random loops whose intensity is given by the
loop measure defined by the semigroup $P_{t}$. More precisely, Propositions
\ref{xx} and \ref{xy} still hold ($\chi$ being defined as a Radon measure with
compact support on $X$) as long as the continuous Green function $G$\ will be
locally trace class. This will apply to exemples related to one dimensional
Brownian motion or to Markov chains on countable spaces.

When points are polar, one needs to be more careful. We will consider only the
case of the two and three dimensional Brownian motion in a bounded domain
killed at the boundary, i.e. associated with the classical energy with
Dirichlet boundary condition. The Green function is not locally trace class
but it is still Hilbert-Schmidt which allows to define renormalized
determinants $\det_{2}$ (Cf \cite{Sim}) and to extend the statement
of\ proposition \ref{xy} to the centered occupation field and the Wick square
$:\phi\overline{\phi}:$\ of the generalized Gaussian Markov field $\phi$.
These three generalized fields are not defined pointwise but have to be
smeared by measures of finite energy $\chi$ such that $\int G^{x,y}%
\chi(dx)\chi(dy)<\infty$. The centered occupation field $\widetilde{l}$ is
defined as follows: Let $A_{t}^{\chi}$ be the additive functional associated
with $\chi$ of finite energy. Then $\left\langle \widetilde{l},\chi
\right\rangle $ is defined as $\lim_{\varepsilon\downarrow0}\int
_{0}^{(T-\varepsilon)^{+}}dA_{t}^{\chi}-\mu_{0}(\int_{0}^{(T-\varepsilon)^{+}%
}dA_{t}^{\chi})$\ which converges in $L^{2}(\mu_{0})$. It is an intrinsic
quantity.. We then have

\begin{proposition}
a) The centered occupation field $\widetilde{\mathcal{L}_{1}}$ and the Wick
square $\frac{1}{2}:\phi\overline{\phi}:$ have the same distribution.

b) $E(e^{-\left\langle \widetilde{\mathcal{L}_{\alpha}},\chi\right\rangle
})=\det_{2}(G_{\chi}G^{-1})^{\alpha}$
\end{proposition}

To justify the use of $\det_{2}$, note that in the finite case $\det
(I+GM_{\chi})=\det_{2}(I+GM_{\chi})e^{\chi(\sigma)}$ where we recall that
$\sigma_{x}^{-1}$ is the capacity of $x$, which vanishes now since $x$ is polar.

In two dimensions, higher Wick powers of $\phi\overline{\phi}$\ are associated
with self intersection local times of the loops.

Let us now consider currents. We will restrict our attention to the two
dimensional Brownian case, $X$ being an open subset of the plane. Currents can
be defined by divergence free vector fields, with compact support. Then
$\int_{l}\omega$ and $\int_{X}(\overline{\phi}\partial_{\omega}\phi
-\phi\partial_{\omega}\overline{\phi})dx$ are well defined square integrable
variables (it can be checked easily in the case of the square by Fourier
series). The distribution of the centered occupation field of the loop process
''twisted'' by the complex exponential $\exp(\sum_{l\in\mathcal{L}_{\alpha}%
}\int_{l}i\omega+\frac{1}{2}\widehat{l}(\left\|  \omega\right\|  ^{2}))$
appears to be the same as the distribution of $:\phi\overline{\phi}%
:$\ ''twisted'' by the complex exponential $\exp(\int_{X}(\overline{\phi
}\partial_{\omega}\phi-\phi\partial_{\omega}\overline{\phi})dx)$ (Cf\cite{NS}).

These points, among others, will be developped in a forthcoming article.

\end{document}